\documentclass{amsart}

\usepackage{amssymb}

\usepackage{amsfonts}
\usepackage{amsmath}

\usepackage{amsthm}
\usepackage{graphicx}
\usepackage{hyperref}

\theoremstyle{plain}

\newtheorem{theorem}{Theorem}
\newtheorem{case}{Case}
\newtheorem{subcase}{Case}
\numberwithin{subcase}{case}

\newtheorem{corollary}{Corollary}

\newtheorem{lemma}{Lemma}

\theoremstyle{definition}

\newtheorem{definition}{Definition}
\newtheorem{remark}{Remark}

\begin{document}

\title{Translation covers among triangular billiard surfaces}

\author{Jason Schmurr}

\address{Department of Mathematics\\
Dalton State College\\
Dalton, GA\\
USA}

\email{jschmurr@daltonstate.edu}

\subjclass[2010]{57M50 (primary); 30F30, 37D50 (secondary)}

\begin{abstract}
We identify all translation covers among triangular billiard
surfaces. Our main tools are the holonomy field of Kenyon and Smillie and a geometric property of triangular billiard surfaces,
which we call the fingerprint of a point, that is preserved under balanced translation
covers.

\end{abstract}

\maketitle

\section{Introduction}\label{section introduction}

An unfolding construction, already described in \cite{Fox Kershner}
and furthered in \cite{Katok Zemlyakov}, associates a flat surface called a \emph{translation
surface} to each rational-angled triangle; we call such a surface a
\emph{triangular billiard surface}. Informally, a compact translation surface is a finite union of polygons in the plane, with pairs of parallel edges identified in such a way as to produce a compact surface. 

Triangular billiard surfaces are highly symmetric examples of translation surfaces; as such, they are of interest in the field -- see Kenyon and Smillie\cite{Kenyon Smillie} and Aurell and Itzykson\cite{Aurell Itzykson}, for example. Structure-preserving maps called \emph{translation covers} between translation surfaces have been used by Vorobets \cite{Vorobets}, Hubert and Schmidt \cite{Hubert Schmidt
2}, Kenyon and Smillie \cite{Kenyon Smillie}, Gutkin and Judge \cite{Gutkin Judge} and others, to gain information about the affine symmetry groups of the surfaces. 

Beyond genus 1, translation covers between triangular billiard surfaces are rare; in fact, our main result is the
following.

\begin{theorem}\label{classification theorem}
Let $f:X \rightarrow Y$ be a nontrivial translation cover of
triangular billiard surfaces, where $X$ has genus greater than 1.
Then each of $X$ and $Y$ is either a right triangular billiards
surface or an isosceles triangular billiard surface, and $f$ is of
degree at most 2.
\end{theorem}

We give an explicit description of which $X$ and $Y$ are related by such covers in Lemma
\ref{right-isosceles lemma}. To prove Theorem \ref{classification
theorem}, we use two main tools: the holonomy field of Kenyon and
Smillie \cite{Kenyon Smillie}, and what we call the
\textit{fingerprint} of a point $P$ on a translation surface, which is an invariant which depends on the surface and a point on the surface and which changes in a natural way under certain translation covers called \emph{balanced covers}. In particular, for a given triangular billiard surface $Y$, we use the holonomy field to narrow the search for possible translation coverings of $Y$ to a finite set of covering surfaces. We then use information (gleaned from fingerprints of points) about the geometric configuration of singular points on these surfaces to identify the actual translation coverings of $Y$\,. 

Our hope is that this strategy of pairing the global information of the holonomy field with the local information of the fingerprint can be applied to other questions involving translation coverings. For example, our main theorem can be seen as establishing a case of a more general conjecture (communicated to us by Eugene Gutkin), that aside from a few trivial cases, most translation coverings between polygonal billiard surfaces are induced from the situation where one of the underlying polygons tiles the other.

\subsection{Outline}

In Sections \ref{subsection billiards construction} and \ref{subsection translation covers}, we review the
construction of triangular billiard surfaces and the basics of translation covers. In Section \ref{subsection the
possible covers}, we prove Lemma \ref{right-isosceles lemma}, which explicitly identifies all possible translation
covers among triangular billiard surfaces. The goal of the remainder of the paper is to show that the list given in Lemma \ref{right-isosceles lemma} is complete -- this is Theorem \ref{classification theorem}. In Section
\ref{section holonomy field} we discuss the holonomy field of Kenyon and Smillie. We offer a new elementary computation of the holonomy field of a given triangular billiard surface, and explain why such surfaces related by a translation cover have the same holonomy field.
In Section \ref{section
fingerprint} we define the fingerprint of a point and prove results
about its behavior under balanced translation covers. 
In \ref{subsection balanced
covers}, we prove the main theorem restricted to balanced covers.
After some combinatorial lemmas in \ref{subsection
combinatorial lemmas}, we complete the proof of the main theorem in \ref{subsection
the main theorem}.

\subsection{Acknowledgments}

The author thanks Thomas A. Schmidt for many
helpful discussions. Thanks are also due to the referee for several useful critiques.

\section{Rational Billiards and Translation Covers}
\subsection{The rational billiards construction}\label{subsection billiards construction}

Let R be a polygonal region whose interior angles are rational
multiples of $\pi$\,. Let $D_{2Q}$ be the dihedral group of order
$2Q$ generated by Euclidean reflections in the sides of $R$\,.
Suppose a particle moves within this region at constant speed and
with initial direction vector $v$, changing directions only when it
reflects off the sides of $R$, with the angle of incidence equaling
the angle of reflection. Every subsequent direction vector for the
particle is of the form $\delta \cdot v$, for some element $\delta
\in D_{2Q}$, where $D_{2Q}$ acts on $\mathbb{R}^2$ via Euclidean
reflections.

The rational billiards construction consists of a flat surface
corresponding to this physical system. Consider the set $D_{2Q}
\cdot R$ of $2Q$ copies of $R$ transformed by the elements of
$D_{2Q}$\,. For each edge $e$ of $R$, we consider the corresponding
element $\rho_e \in D_{2Q}$ which represents reflection across
$e$\,. For each $\delta \in D_{2Q}$, we glue $\rho_e \delta \cdot R$
and $\delta \cdot R$ together along their copies of $e$\,. The
result is a closed Riemann surface with flat structure induced by
the tiling by $2Q$ copies of $R$\,. This construction is described
by Fox and Kershner in \cite{Fox Kershner}.

In fact this surface is an example of a compact \emph{translation surface}. A compact translation surface can be defined as the result of gluing together a finite set of polygons in the plane along parallel edges in such a way that the result is a compact surface (see, e.g., \cite{Calta Smillie}). 

Equivalently, a translation surface can be defined as a real two-manifold. Recall that given two coordinate maps $\phi_1:U_1\rightarrow \mathbb{R}^2$ and $\phi_2:U_2\rightarrow \mathbb{R}^2$ defining homeomorphisms from open sets $U_1$ and $U_2$ of a manifold $M$ into $\mathbb{R}^2$, the map $\phi_2\circ\phi_1^{-1}:\phi_1(U_1\cap U_2)\rightarrow \phi_2(U_1\cap U_2)$ is called a \emph{transition function}. A conical singularity $P$ on a flat surface is a point such that, in the flat metric induced by the coordinate maps, the total angle (``cone angle") about $P$ is not equal to $2\pi$\,. 

\begin{definition}
Let $X$ be a flat surface with conical singularities. Let
$\tilde{X}$ be the flat surface obtained by puncturing all
singularities of $X$\,. If all transition functions of $\tilde{X}$
are translations, then $X$ is a \emph{translation surface}.
\end{definition}

On a translation surface, the cone angles of conical singularities are always integer multiples of $2\pi$\,. See \cite{Zorich} for an introduction to flat surfaces.

We focus on billiards in rational triangles. Let $\left(a_1,a_2,a_3\right)$ be a triple of positive integers. We fix the notation
$T(a_1,a_2,a_3)$ to refer to a triangle with internal angles
$\dfrac{a_1\pi}{Q}$, $\dfrac{a_2\pi}{Q}$, and $\dfrac{a_3\pi}{Q}$,
where $Q:=a_1+a_2+a_3$ and $\gcd(a_1,a_2,a_3)=1$\,. We use the
notation $X(a_1,a_2,a_3)$ to refer to the translation surface
arising from billiards in $T(a_1,a_2,a_3)$ via the Fox-Kershner
construction. We call such a surface a \emph{triangular billiards
surface}. If the triangle is isosceles or right, we call the
corresponding surface an \emph{isosceles triangular billiards
surface} or a \emph{right triangular billiard surface}.

\begin{definition}
Note that the Fox-Kershner construction gives a natural ``tiling by
flips'' of the surface by copies of $T$\,. A \emph{billiards
triangulation} is a triangulation $\tau$ of $X$ whose triangles are
the various elements of $D_{2Q} \cdot T$ described above.
\end{definition}

\begin{remark}\label{remark on AI combinatorics}
Some data can be gained about $\tau$ via simple combinatorics.
Letting $T:=T(a_1,a_2,a_3)$, label the vertices of $T$ as $v_1$,
$v_2$, and $v_3$, where $v_i$ corresponds to $a_i$\,. It is not hard
to check that the total number of triangles in $\tau$ is $2Q$, that
the number of vertices of $\tau$ corresponding to $v_i$ is
$\gcd(a_i,Q)$, and that each member of this set has a cone angle of
$\left(\dfrac{a_{i}}{\gcd(a_i,Q)}\right)2\pi$\,. A good reference
for these matters is \cite{Aurell Itzykson}.
\end{remark}

\begin{definition}
Let $v_1,v_2,v_3$ be the vertices of  $T(a_1,a_2,a_3)$, and let
\newline $\pi_X:X(a_1,a_2,a_3) \rightarrow T(a_1,a_2,a_3)$ be the
standard projection. A \emph{vertex class} of $X(a_1,a_2,a_3)$ is
any of the three sets $\pi_X^{-1}(v_1)$, $\pi_X^{-1}(v_2)$, or
$\pi_X^{-1}(v_3)$\,. Note that for a given vertex class, either all
the elements are singular or all are nonsingular; hence we call a
vertex class \emph{singular} if its elements are singularities and
\emph{nonsingular} if its elements are nonsingular.
\end{definition}

Clearly, a vertex class $\pi_X^{-1}(v_i)$ is nonsingular if and only
if $a_i$ divides $Q$\,. Furthermore, the sum of the cone angles of
the elements of $\pi_X^{-1}(v_i)$ is $2a_i\pi$\,. 

\subsection{Translation covers}\label{subsection translation covers}
The natural map between translation surfaces is one which respects the
translation structure:

\begin{definition}
A \emph{translation cover} is a holomorphic (possibly ramified)
cover of translation surfaces $f:X \rightarrow Y$ such that, for
each pair of coordinate maps $\phi_X$ and $\phi_Y$ on $X$ and $Y$,
respectively, the map $\phi_Y \circ f \circ \phi_X^{-1}$ is a
translation when $\phi_X$ and $\phi_Y$ are restricted to open sets not containing singular
points. We say that $f$ is \emph{balanced} if $f$ does not map singular points to nonsingular points.
\end{definition}

The term ``balanced" cover is due to Gutkin \cite{Gutkin}. Balanced translation covers $f: X \rightarrow Y$ of translation
surfaces are of particular interest because they imply an especially strong
relationship between the affine symmetry groups of $X$ and $Y$; in
particular, these groups must share a subgroup which is of finite index in each (see  \cite{Gutkin Judge} and \cite{Vorobets}).

\begin{definition}
We say that $X$ and $Y$ are \emph{translation equivalent} if there
exists a degree 1 translation cover $f:X \rightarrow Y$\,.
\end{definition}

\begin{figure}[ht]
\includegraphics[scale=0.4]{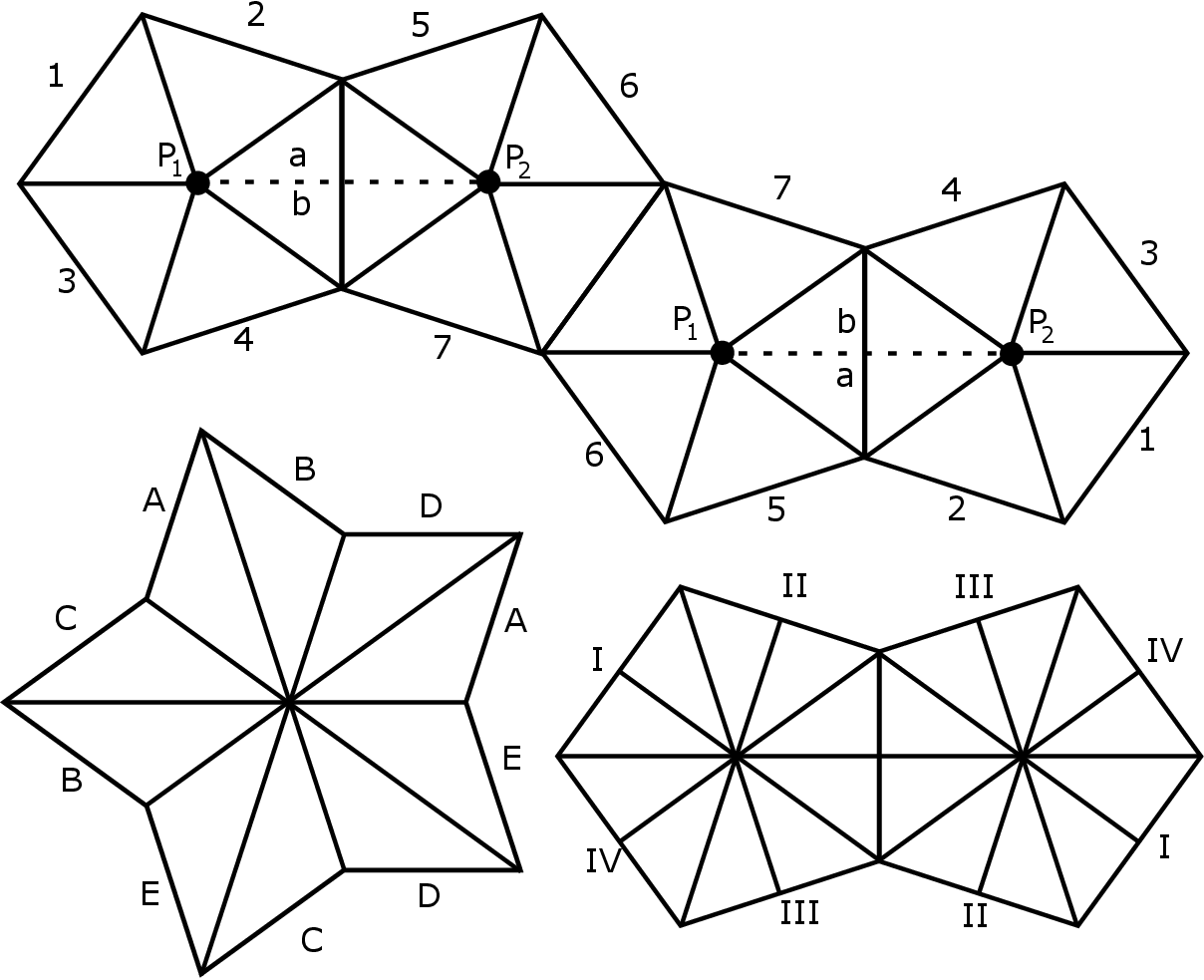}
\caption{Above, the surface $X(4,3,3)$\,. Below, $X(3,1,1)$ (left) and $X(5,3,2)$ (right)\,. }
\label{433532311}
\end{figure}

The surface $X(4,3,3)$ is pictured at the top of Figure \ref{433532311}. The emphasized points labeled $P_1$ and $P_2$ form the singular vertex class corresponding to the vertex angle of $2\pi/5$ in $T(4,3,3)$\,. Each of these points has a cone angle of $4\pi$, because 10 copies of $T(4,3,3)$ are developed about each point. To help visualize this, we follow Hubert and Schmidt \cite{Hubert Schmidt 3} and use the dotted lines to represent slits with gluings indicated by labels $a$ and $b$\,. The identification of parallel edges as indicated by the labeling scheme in the figure leads to two more singular points on $X(4,3,3)$, each of cone angle $6\pi$\,. Also in Figure \ref{433532311} we depict $X(3,1,1)$ and $X(5,3,2)$\,. These two surfaces are translation equivalent, and there exist degree two translation covers from $X(4,3,3)$ to each of them.

The following lemma demonstrates how we will use Remark \ref{remark on AI combinatorics}
to analyze translation covers.

\begin{lemma}\label{lemma cover combinatorics}
Suppose $f:X(a_1,a_2,a_3) \rightarrow X(b_1,b_2,b_3)$ is a
translation cover of triangular billiard surfaces. Let
$\pi_X:X(a_1,a_2,a_3) \rightarrow T(a_1,a_2,a_3)$ and
$\pi_Y:X(b_1,b_2,b_3) \rightarrow T(b_1,b_2,b_3)$ be the canonical
projections to triangles with vertices $v_1,v_2,v_3$ and
$w_1,w_2,w_3$ respectively. Suppose that $P \in \pi_Y^{-1}(w_i)$,
$P' \in \pi_X^{-1}(v_j)$, and $f(P')=P$ with a ramification index of
$m$ at $P'$\,. Then

$\dfrac{m b_i}{\gcd(b_i,b_1+b_2+b_3)}=
\dfrac{a_j}{\gcd(a_j,a_1+a_2+a_3)}$\,.
\end{lemma}

\begin{proof}
The cone angle at $P'$ is $m$ times the cone angle at $P$\,.
Therefore the result follows from Remark \ref{remark on AI
combinatorics}.
\end{proof}

Of course, the translation structure of $X(a_1,a_2,a_3)$ depends on
the chosen area and direction of $T(a_1,a_2,a_3)$\,. A translation
surface $X$ can represented as a pair $(S,\omega)$, where $S$ is a
Riemann surface and $\omega$ is a holomorphic 1-form on $S$ which
induces the translation structure of $X$\,. Using this language,
suppose that $(S,\omega)$ is a triangular billiard surface arising
from billiards in some $T(a_1,a_2,a_3)$, and that $\alpha$ is a
nonzero complex number. The notation $X(a_1,a_2,a_3)$ does not
distinguish the pairs $(S,\omega)$ and $(S,\alpha\omega)$\,. The
following lemma shows that this ambiguity will not affect our
classification of translation covers.

\begin{lemma}\label{singular surface self-cover lemma}
Suppose that $(S,\omega)$ is a triangular billiard surface of genus
greater than one, and let $\alpha \in \mathbb{C} \backslash
\{0\}$\,. Then any translation cover $f:(S,\omega) \rightarrow
(S,\alpha\omega)$ is of degree 1.
\end{lemma}

\begin{proof}
This is a simple application of the Riemann-Hurwitz formula. Let
$(S,\omega)$ have genus $g$, and let $\deg f =n$\,. The 1-form
$\omega$ which gives $(S,\omega)$ its translation structure has
$2g-2$ zeros (counting multiplicities). Clearly $\alpha\omega$ has
the same zeros as $\omega$\,. The Riemann-Hurwitz formula then gives
us that
\begin{equation}\label{RH calculation}
g=n(g-1)+1+\dfrac{R}{2},
\end{equation}
where $R$ is the total ramification number of $f$\,. Since $R \geq
0$, Equation \eqref{RH calculation} is only satisfied if $n=1$\,.
\end{proof}

\subsection{The Possible Translation Covers}\label{subsection the possible covers}
Any isosceles triangle is naturally ``tiled by flips'' by a right
triangle. The following lemma demonstrates how to use this tiling to
create nontrivial translation covers in the category of triangular
billiard surfaces. In fact, our main theorem is that the covers of
Lemma \ref{right-isosceles lemma} are the only nontrivial
translation covers among triangular billiard surfaces.

\begin{lemma}\label{right-isosceles lemma}
Let $a$ and $b$ be relatively prime positive integers, not both equal to one. The right triangular billiard surface $Y:=X(a_1+a_2,a_1,a_2)$ is related to two isosceles triangular
billiard surfaces $X_1$ and $X_2$ via translation covers $f_i:X_i \rightarrow Y$\,. For each $i$, if $a_i$ is odd then $X_i=X(2a_j,a_i,a_i)$ and $f_i$ has degree 2. If $i$ is even then $X_i=X\left(a_j,\dfrac{a_i}{2},\dfrac{a_i}{2}\right)$ and $f_i$ has degree 1.
\end{lemma}

\begin{proof}
It suffices to prove the result for $X_1$ and $f_1$\,. Write
$Q:=2a_1+2a_2$\,. We reflect the triangle $T=T(a_1+a_2,a_1,a_2)$
across the edge connecting the vertices corresponding to $a_2$ and $a_1+a_2$, to
obtain its mirror image $T'$\,. By joining $T$ and $T'$ along the
edge of reflection we create an isosceles triangle $\tilde{T}$ which
can be written as either $T(2a_2,a_1,a_1)$ (if $a_1$ is odd)  or
$T(a_2,\dfrac{a_1}{2},\dfrac{a_1}{2})$ (if $a_1$ is even). Note that
since $(a_1+a_2,a_1,a_2)$ must be a reduced triple, $a_1$ and $a_2$
cannot both be even. It also follows that $\gcd(a_i,Q) \leq
\gcd(2a_i,Q) = 2$\,.

Suppose $a_1$ is even. Consider the translation surface $S$ (with
boundary) obtained by developing $T$ around its vertex corresponding to $a_2$\,. Since
$a_2$ is odd we have $\gcd(a_2,Q)=1$, so $S$ is tiled (by
reflection) by $2Q$ copies of $T$, and hence after appropriate
identifications along the boundary we will have
$X(a_1+a_2,a_1,a_2)$\,. Let $\tilde{S}$ be the surface obtained by
developing $\tilde{T}$ around the corresponding vertex; it is tiled
via reflection by $Q$ copies of $\tilde{T}$, so appropriate boundary
identifications will yield $Y_1$\,. Because $\tilde{T}$ is tiled via
reflection by two copies of $T$, it follows that $S$ and $\tilde{S}$
are translation equivalent. Finally, note that the boundary
identifications are the same for $S$ and $\tilde{S}$\,. Therefore
$Y$ and $X_1$ are translation equivalent.

Now suppose that $a_1$ is odd and $a_2$ is even. We then have
$\tilde{T}=T(2a_2,a_1,a_1)$\,. Since $\gcd(2a_2,Q)=2$, we again have
that $\tilde{S}$ is tiled by $Q$ copies of $\tilde{T}$\,. Since
$a_2$ is even, $\gcd(a_2,Q)=2$, implying that $S$ is tiled by $Q$
copies of $T$\,. Thus if $a_2$ is even then there exists a degree
two cover $f: \tilde{S} \rightarrow S$\,, ramified over a single
point. Furthermore, in this case $X_1$ and $Y$ are obtained by
identifying appropriate edges of \emph{two} copies of $\tilde{S}$
and $S$, respectively. It follows that if $a_2$ is even then there
exists a ramified degree two cover $f:X_1 \rightarrow Y$\,.

Finally, suppose that $a_1$ and $a_2$ are both odd. We have
$\tilde{T}=T(2a_2,a_1,a_1)$, $\gcd(2a_2,Q)=2$, and
$\gcd(a_2,Q)=1$\,. In this case we have that $S$ and $\tilde{S}$ are
translation equivalent surfaces; however, $X_1$ is obtained from two
copies of $\tilde{S}$ whereas $Y$ is obtained from a single copy of
$S$\,. Thus again we have that $X_1$ is a double cover of $Y$, this
time unramified.
\end{proof}

\begin{remark}Note that in addition to relating right and isosceles triangles,
Lemma \ref{right-isosceles lemma} also gives a way to construct
covers between isosceles triangular billiard surfaces. In the
language of Lemma \ref{right-isosceles lemma}, if $a_2$ is even,
then $f_2^{-1} \circ f_1$ is a degree two translation cover of $X_2$
by $X_1$\,.
\end{remark}

\begin{remark}\label{remark on genus 1}
If we allow $a_1=a_2=1$  in the statement of Lemma
\ref{right-isosceles lemma}, then we arrive at
$Y=X_1=X_2=X(1,1,2)$\,. This is because $T(1,1,2)$ is the unique
right isosceles triangle.
\end{remark}

Because the location of singularities is such a major tool in
analyzing translation surfaces, it is worth identifying the
triangular billiard surfaces which have no singularities. As
detailed in \cite{Aurell Itzykson}, there are only three of these
surfaces: $X(1,1,2)$, $X(1,2,3)$, and $X(1, 1, 1)$\,. These are also
the only three triangular billiard surfaces of genus 1; furthermore
$X(1,2,3)$ and $X(1,1,1)$ are actually translation equivalent. Each
of these surfaces admits balanced translation covers of itself by
itself of arbitrarily high degree; this fact is related to the fact
that $T (1,1,2)$, $T (1,2,3)$, and $T (1,1,1)$ are the only
Euclidean triangles which tile the Euclidean plane by flips. Note
that any such cover must be unramified, since about ramification points flat ramified covers
are locally of the form $z \mapsto z^{1/n}$ for some $n>1$, implying
that the cone angle of the ramification point is greater than
$2\pi$; hence ramification points are singular\,.

\section{The holonomy field}\label{section holonomy field}

In order to prove that Lemma \ref{right-isosceles lemma} provides a complete list of the translation covers between triangular billiard surfaces, we require some information about the holonomy of a triangular billiard surface.

\begin{definition}\label{absolute holonomy definition}
The rational \emph{absolute holonomy} of a translation surface $X$
is the image of the map $\text{hol}:H_1(X;\mathbb{Q}) \rightarrow
\mathbb{C}$ defined by $\text{hol}:\sigma \mapsto
\int_{\sigma}\omega$, where the 1-form $\omega$ is locally the
differential $dz$ in each chart not containing a singular point.
\end{definition}

The following definition is due to Kenyon and Smillie \cite{Kenyon Smillie}.

\begin{definition}\label{holonomy field definition}
The \emph{holonomy field} of a translation surface $X$, denoted
$k_X$, is the smallest field $k_X$ such that the rational absolute holonomy
of $X$ is contained in a two-dimensional vector space over $k_X$\,.
\end{definition}

\subsection{Calculation of the holonomy field} \label{subsection holonomy calculation}

Kenyon and Smillie \cite{Kenyon Smillie} calculate the holonomy field of $X(a_1,a_2,a_3)$ to be $k_X=\mathbb{Q}\left(\cos \left(2\pi/{Q}\right)\right)$\,. We offer a more elementary proof of this result. 

\begin{remark}\label{remark:Chebyshev}
We let $U_n$ denote the $n^{th}$ Chebyshev polynomial of the second kind. We will use the following properties of Chebyshev polynomials.
\begin{enumerate}
  \item $\frac{\sin((n+1)\theta)}{\sin\theta}=U_{n}(\cos\theta)$
  \item If $n$ is even, then $U_{n}$ is an even polynomial of degree $n$\,. If $n$ is odd, then $U_{n}$ is an odd polynomial of degree $n$\,.
\end{enumerate}
\end{remark}

\begin{remark}\label{remark:cyclotomic equality}
Let $\phi$ be the Euler totient function. It is well known that, for
any positive integer $Q$, the degree of the number field
$\mathbb{Q}(\cos(2\pi/Q))$ is equal to $\frac{1}{2}\phi(Q)$\,. Note that
if $Q$ is odd, then $\phi(Q)=\phi(2Q)$\,. It follows that, when $Q$ is
odd, we will have
$\mathbb{Q}\left(\cos \left(2\pi/{Q}\right)\right)=\mathbb{Q}\left(\cos \left(\pi/{Q}\right)\right)$\,.
\end{remark}

The following is Proposition 2.5 of \cite{Calta Smillie}.
\begin{lemma}\label{thm:holonomy in Lambda}(Calta-Smillie)
If a translation surface $X$ is obtained by identifying the edges of
polygons in the plane by maps which are restrictions of
translations, and if all the vertices of these polygons lie in a
subgroup $\Lambda \subset \mathbb{R}^{2}$, then the holonomy of $S$
is contained in $\Lambda$\,.
\end{lemma}

\begin{figure}[ht]
\includegraphics[scale=0.4]{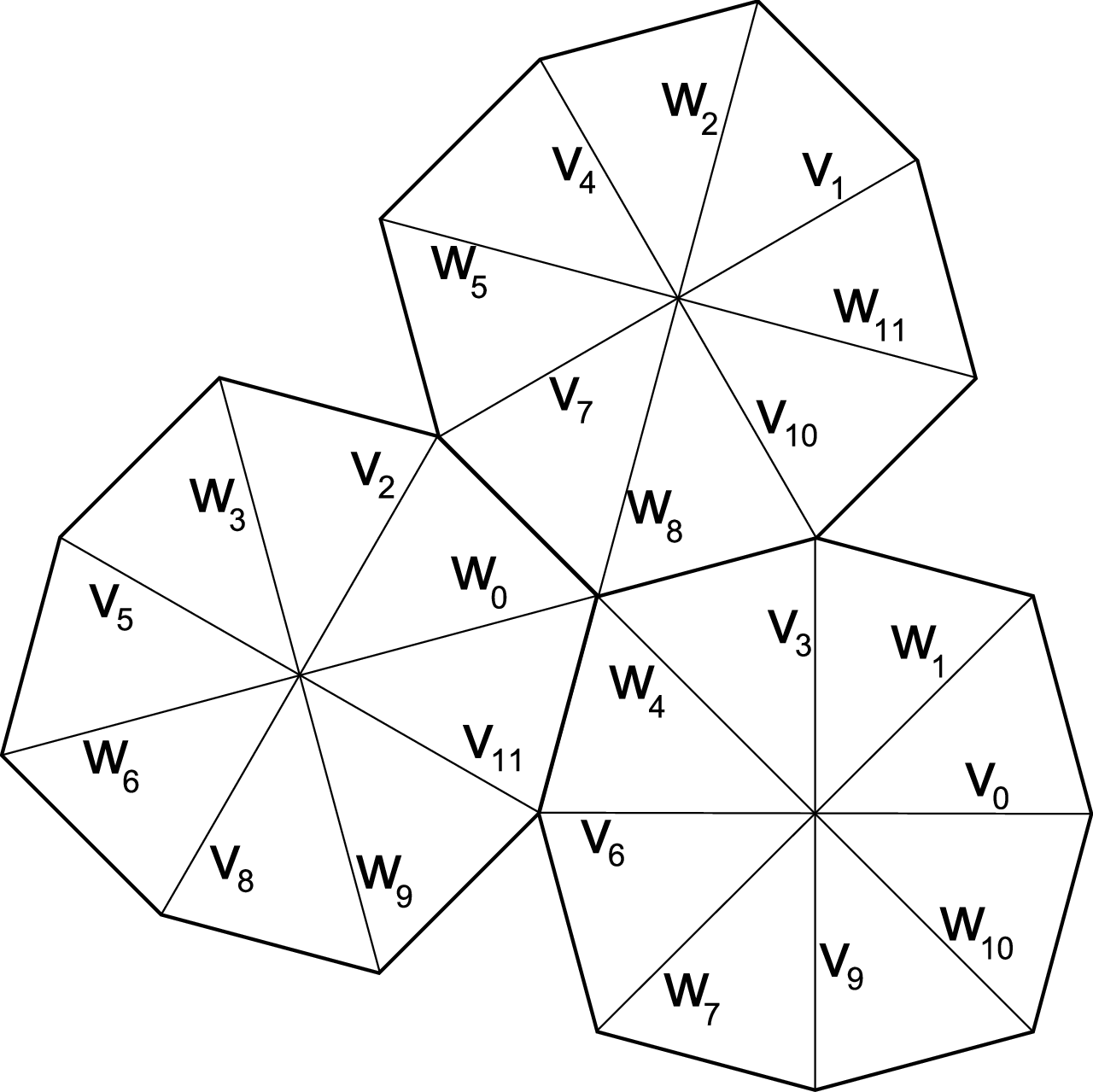}
\caption{The sets $\{v_n\}$ and $\{w_n\}$ for $X(3,4,5)$, with $a_1=3$\,.}\label{figure:X(3,4,5)_algper}
\end{figure}

\begin{lemma}\label{holonomy field calculation lemma}

The holonomy field of $X=X(a_1,a_2,a_3)$ is $k_X=\mathbb{Q}\left(\cos \left(2\pi/{Q}\right)\right)$\,.
 
\end{lemma}

\begin{proof}
Let $k_X$ denote the holonomy field of $X$\,.
Let $\alpha=\frac{\pi}{Q}$\,. Let $T=T(a_1,a_2,a_3)$\,. Since $\gcd(a_1,a_2,a_3)=1$, we can and do assume that $a_1$ is odd. Label the vertices of $T$ corresponding to the angles
$a_{1}\alpha$, $a_{2}\alpha$, and $a_{3}\alpha$ as $P_{1}$, $P_{2}$,
and $P_{3}$\,. We scale and rotate $T$ so that the
$\overline{P_{1}P_{2}}$ side has edge vector $v=(1,0)$, and so that
the $\overline{P_{1}P_{3}}$ side has edge vector
$w=\left(t\cos(a_{1}\alpha),t\sin(a_{1}\alpha)\right)$, where by the Law of Sines we have $t=\dfrac{\sin(a_{2}\alpha)}{\sin(a_{3}\alpha)}$\,. The dihedral group
$D$ generated by reflections in the sides of $T$ acts on the
set $D\cdot T$ of $2Q$ distinct oriented triangles arising from
billiards in $T$\,. We can construct $X$ from this set by
identifying the appropriate edges of the elements of $D\cdot
T$\,. We may also view $D$ as acting on the edge vectors of
$T$\,. Let $v_{n}=\left(\cos(2n\alpha),\sin(2n\alpha)\right)$ and
$w_{n}=(t\cos((2n+1)\alpha),t\sin((2n+1)\alpha))$\,. With this
notation, we see that $D \cdot v$ is the set
$\{v=v_{0},v_{1},...,v_{Q-1}\}$\,. Recalling that $a_{1}$ is odd, we also see that $D \cdot w$ is the set $\{w_{0},w_{1},...,w_{Q-1}\}$\,. Note that $w=w_{(a_1-1)/2}$\,.

Let $L=\mathbb{Q}(\cos (2\alpha))$\,. We will show that all the $v_{n}$
and $w_{n}$ are $L$-linear combinations of $v_{0}$ and $v_{1}$, and
that furthermore $L$ is the smallest such field.

To see that $k_X$ contains $L$, we note that $v_2$ is the result of rotating $v_1$ by an angle of $2\alpha$\,. Consider a vector space over $k_X$ with ordered basis $\left\{v_0,v_1\right\}$\,. In this basis, counterclockwise rotation by $2\alpha$ is a linear transformation $R$ with matrix $\left(
                                                                                                                               \begin{array}{cc}
                                                                                                                                 0 & -1 \\
                                                                                                                                 1 & 2\cos(2\alpha) \\
                                                                                                                               \end{array}
                                                                                                                             \right)$\,. Hence, $\cos(2\alpha)$ is in $k_X$, and we have $L \subseteq k_X$\,.

Since each $v_i$ and each $w_i$ is the result of repeated application of $R$ to $v_0$ or $w_0$, it remains to be shown that $w_0$ is an $L$-linear combination of $v_0$ and $v_1$\,. Then Lemma \ref{thm:holonomy in Lambda} gives that $k_X \subseteq L$, completing the proof.

Let $l$ and $l'$ be the real numbers such that
$lv_{0}+l'v_{1}=w_{0}$\,. Since $v_{0}$ and $v_{1}$ are reflections of
each other across the line generated by $w_{0}$, we see that
$v_{0}+v_{1}$ is a real multiple of $w_{0}$\,. Hence $l'=l$\,.

Projecting $v_{0}$ and $v_{1}$ onto $w_{0}$, we see that

\begin{equation}
l=\frac{||w_{0}||}{||v_{0}+v_{1}||} =\frac{t}{2\cos\alpha}
=\frac{\sin(a_{2}\alpha)\sin\alpha}{\sin(a_{3}\alpha)\sin(2\alpha)}
=\frac{\sin(a_{2}\alpha)}{\sin\alpha}\frac{\sin\alpha}{\sin(a_{3}\alpha)}\frac{\sin\alpha}{\sin(2\alpha)}.
\end{equation}

Applying Remark \ref{remark:Chebyshev} to the last expression, we get

\begin{equation}
l=\frac{U_{a_{2}-1}(\cos\alpha)}{U_{a_{3}-1}(\cos\alpha)U_{1}(\cos\alpha)}.
\end{equation}

If $Q$ is even, we have that $(a_{2}-1)$ and $(a_{3}-1)$ have
opposite parity, and thus by our Remark \ref{remark:Chebyshev},
$\frac{U_{a_{2}-1}(\cos\alpha)}{U_{a_{3}-1}(\cos\alpha)U_{1}(\cos\alpha)}$
is a rational function in $\cos^{2}\alpha$\,. Therefore  $l\in
\mathbb{Q}(\cos^{2}\alpha)=L$\,. If $Q$ is
odd, then already by Remark \ref{remark:cyclotomic equality},
 $\mathbb{Q}(\cos\alpha)=L$, and since $\frac{U_{a_{2}-1}(\cos\alpha)}{U_{a_{3}-1}(\cos\alpha)U_{1}(\cos\alpha)}$
is a rational function in $\cos\alpha$, we again have that $l \in
L$\,.

Hence
$\text{span}_{L}\{v_{0},v_{1}\}$=$\Lambda$\,. Theorem \ref{thm:holonomy in Lambda} says that
$\Lambda$ contains the absolute holonomy of $S$\,. So $L$ contains $k_X$, which completes the proof.
\end{proof}

\begin{corollary}\label{holonomy field equality corollary}
If $f:X \rightarrow Y$ is a translation cover of triangular billiard surfaces, then $k(X)=k(Y)=\mathbb{Q}(\cos(2\pi/Q))$\,.
\end{corollary}
\begin{proof}
Since $f$ is a translation cover, the billiards triangulation $\tau$ of $Y$ by $T_Y$ lifts to a triangulation of $X$ by $T_Y$, in which each edge is is mapped by $f$ to an edge of $\tau$\,. Hence we may take the subgroup $\Lambda$ referred to in Lemma \ref{thm:holonomy in Lambda} to be the same group for both $X$ and $Y$\,. Then our calculation in Lemma \ref{holonomy field calculation lemma} will be the same for $X$ and $Y$\,.
\end{proof}

\subsection{Some elementary number theory}
Since $k_{X}=\mathbb{Q}(\cos(2\pi/Q))$, we would like to know when distinct values of $Q$ yield the same field $k_X$\,.

Letting $\zeta_Q$ denote a primitive $Q^{th}$ root of unity, we have that $\mathbb{Q}(\cos(2\pi/Q))$ is equal to $\mathbb{Q}(\zeta_{Q}+\zeta_{Q}^{-1})$, which is a degree two
subfield of the cyclotomic field $\mathbb{Q}(\zeta_{Q})$, since it
is the maximal subfield fixed by complex conjugation. In light of
this, we list some classical results about these two fields as
recorded in Washington's text\cite{Washington}.

\begin{lemma}\label{cyclotomic 2 lemma}
If $Q$ is odd then $\mathbb{Q}(\zeta_{Q})=\mathbb{Q}(\zeta_{2Q})$\,.
\end{lemma}

\begin{lemma}\label{Wash 2.3}(Prop 2.3 in \cite{Washington}) Assume that
$Q \not\equiv 2 \text{ mod }4$\,. A prime $p$ ramifies in
$\mathbb{Q}(\zeta_{Q})$ if and only if $p|Q$\,.
\end{lemma}

\begin{lemma}\label{Wash 2.15}(Prop 2.15 in \cite{Washington}) Let $p$ be a prime,
and assume that $Q \not\equiv 2 \text{ mod }4$\,. If $Q=p^{m}$ then
$\mathbb{Q}(\zeta_{Q})/\mathbb{Q}(\zeta_{Q}+\zeta_{Q}^{-1})$ is
ramified only at the prime above $p$ and at the archimedean primes.
If $Q$ is not a prime power, then
$\mathbb{Q}(\zeta_{Q})/\mathbb{Q}(\zeta_{Q}+\zeta_{Q}^{-1})$ is
unramified except at the archimedean primes.
\end{lemma}

\begin{remark}\label{Washington fix}
Washington's proofs of Lemmas \ref{Wash 2.3} and \ref{Wash 2.15}
make clear that the results carry through to the case $Q \equiv 2
\text{ mod } 4$ except that in that case, the prime 2 does not
ramify in $\mathbb{Q}(\zeta_{Q})$\,.
\end{remark}

For a triangular billiard surface $X=X(a_1,a_2,a_3)$, it is
tempting to define a ``$Q$-value'' for the surface by
$Q_{X}:=a_1+a_2+a_3$\,. Unfortunately this notion is not quite
well-defined up to translation equivalence; as demonstrated in Lemma
\ref{right-isosceles lemma}, the distinct triangles $T(a,a,b)$ and
$T(2a,b,2a+b)$ unfold to translation equivalent translation surfaces
if (and only if) $b$ is odd. However, the following lemma and its
corollary show that this notion is well-defined up to a factor of
$2$\,.

\begin{lemma}\label{maximal real subfield lemma}
Distinct cyclotomic fields have distinct maximal totally real
subfields.
\end{lemma}

\begin{proof}
This is an exercise in elementary algebraic number theory, and is
presumably well known. Let $k$ be the maximal totally real subfield of the
cyclotomic fields $\mathbb{Q}(\zeta_m)$ and $\mathbb{Q}(\zeta_n)$
for positive integers $m,n>2$\,. Let $p$ be an odd prime dividing
$m$\,. By Lemma \ref{Wash 2.3}, $p$ ramifies in
$\mathbb{Q}(\zeta_m)$\,. If $m$ is a power of $p$, then $p$ is
totally ramified in $\mathbb{Q}(\zeta_m)$\,. Since
$\mathbb{Q}\subset k \subset \mathbb{Q}(\zeta_m)$, if $m$ is a power of $p$ then
$p$ must ramify in $k$\,. If $m$ is not a power of $p$\,, then Lemma
\ref{Wash 2.15} tells us that the extension $\mathbb{Q}(\zeta_m)/k$
is not ramified at the prime above $p$ ; thus again $p$ must ramify
in $k$\,. But also $\mathbb{Q} \subseteq k \subset \mathbb{Q}(\zeta_n)$, so $p$
must ramify in $\mathbb{Q}(\zeta_n)$\,. By Lemma \ref{Wash 2.3},
this implies that $p$  divides $n$\,. Therefore $m$ and $n$ have the
same odd prime divisors; furthermore, by Remark \ref{Washington
fix}, these arguments extend to show that either 4 divides both $m$
and $n$ or it divides neither.

The degrees of $\mathbb{Q}(\zeta_m)$ and $\mathbb{Q}(\zeta_n)$ as
field extensions of $\mathbb{Q}$ are $\phi(m)$ and $\phi(n)$
respectively, where $\phi$ is the Euler totient function. Since
$\mathbb{Q}(\zeta_m)$ and $\mathbb{Q}(\zeta_n)$ are each degree $2$
extensions of $k$, we have that $\phi(m)=\phi(n)$\,.

First suppose that $m$ and $n$ are congruent modulo $2$\,. Let
$m=\Pi p_{i}^{e_{i}}$ and $n=\Pi p_{i}^{f_{i}}$ be the prime
factorizations of $m$ and $n$\,. Then we have

\begin{equation}\label{totient quotient equation}
1=\dfrac{\phi(m)}{\phi(n)}=\dfrac{\prod
(p_{i}-1)p_{i}^{e_{i}-1}}{\prod (p_{i}-1)p_{i}^{f_{i}-1}}= \prod
p_{i}^{e_{i}-f_{i}}.
\end{equation}

Therefore $e_{i}=f_{i}$ for each $i$, and $m=n$\,. Hence
in this case $\mathbb{Q}(\zeta_m)=\mathbb{Q}(\zeta_n)$\,.

If $m$ and $n$ are not congruent modulo $2$, then we may assume that
$m$ is odd and $n$ is congruent to $2$ modulo $4$\,. Since
$\phi(m)=\phi(2m)$ when $m$ is odd, we can repeat the calculation
\eqref{totient quotient equation} with $2m$ and $n$, and get that
$2m=n$\,. But it is well known that for any odd $m$,
$\mathbb{Q}(\zeta_m)=\mathbb{Q}(\zeta_{2m})$\,. Therefore in fact
$k$ is the maximal totally real subfield of only one cyclotomic
field.
\end{proof}

\begin{corollary}\label{holonomy determines Q corollary}
Suppose that $X_1=X(a_1,a_2,a_3)$ and $X_2=X(b_1,b_2,b_3)$ are related by a translation cover, and that $b_1+b_2+b_3< a_1+a_2+a_3$\,. Then
$b_1+b_2+b_3$ is odd, and $a_1+a_2+a_3=2(b_1+b_2+b_3)$\,.
\end{corollary}

\begin{proof}
 By Corollary \ref{holonomy field equality corollary}, $X_1$ and $X_2$ have the same holonomy
field $k$\,. Write $Q_{X_1}=a_1+a_2+a_3$ and $Q_{X_2}=b_1+b_2+b_3$\,.  Then
 $k$ is the maximal
totally real subfield of $\mathbb{Q}(\zeta_{Q_{X_1}})$ and of
$\mathbb{Q}(\zeta_{Q_{X_2}})$\,. Hence by Lemma 8, $\mathbb{Q}(\zeta_{Q_{X_1}})=\mathbb{Q}(\zeta_{Q_{X_2}})$\,. The result then follows directly from the proof of
Lemma \ref{maximal real subfield lemma}.
\end{proof}

\section{The Fingerprint}\label{section fingerprint}

We have seen in Section \ref{section holonomy field} that we can use the holonomy field to reduce the list of possible translation covers of a given triangular billiard surface $Y$ to a finite set: those surfaces $X$ such that $Q_X$ and $Q_Y$ differ by at most factor of two. To decide which of these surfaces actually are related by translation covers, we will study the geometric configuration of singularities on a billiard surface by examining the shortest paths from points on the surface to singularities.

Consider a point $P$ on a translation surface $X$, along with the
set of $S$ all shortest geodesic segments on $X$ which connect $P$
to a singularity. Let $s_1$ and $s_2$ be two of these segments. We
say that $s_1$ and $s_2$ are \emph{adjacent} if $s_1$ can be
rotated continuously about $P$ onto $s_2$ without first coinciding
with any other elements of $S$\,.

\begin{definition}
A \emph{fingerprint} of a point $P\in\tau$ is the data
$\{\{\theta_{i}\},\phi,d\}$, where $\{\theta_{i}\}$ contains the
distinct angle measures separating adjacent pairs of shortest
geodesic segments connecting $P$ to singularities, $\phi$ is the
total cone angle at $P$, and $d$ is the length of each of the
shortest geodesic segments. We shall say that $P$ has a \emph{Type
I fingerprint} if $\{\theta_{i}\}$ has one element, and that $P$ has
a \emph{Type II fingerprint} if $\{\theta_{i}\}$ has two elements.
We call $\{\theta_{i}\}$ the \emph{angle set} of a fingerprint.
\end{definition}

\begin{figure}[ht]
\includegraphics[scale=.6]{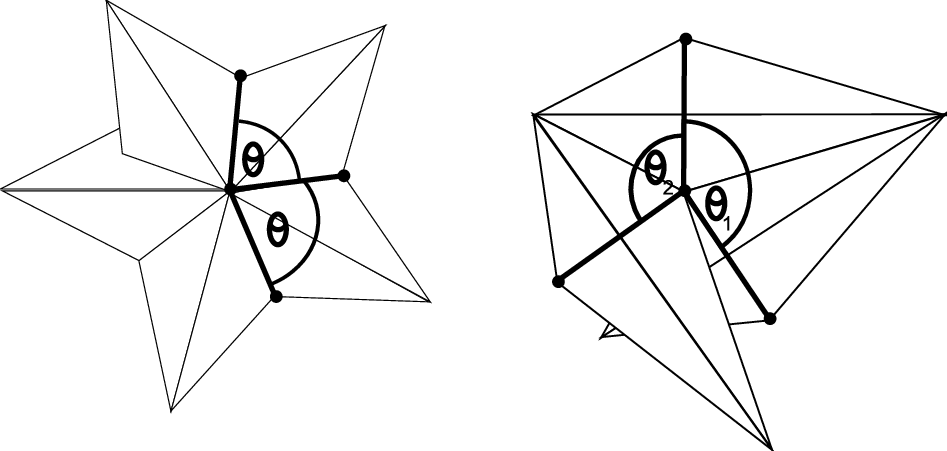}
\caption{Parts of a Type I fingerprint (left) and a Type II
fingerprint (right).}\label{Type I and II figure}
\end{figure}

\begin{figure}[ht]
\includegraphics[scale=.25]{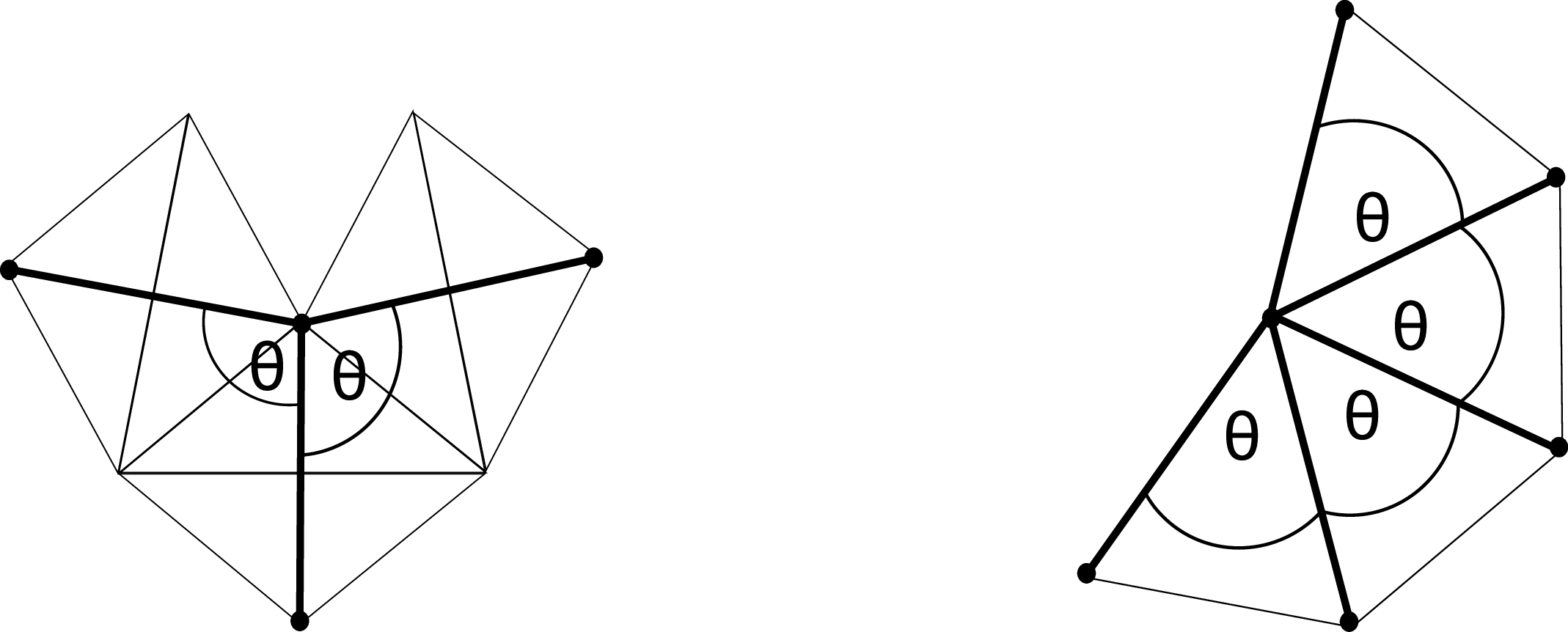}
\caption{Type I fingerprints arising from isosceles triangles.}\label{isosceles type I figure}
\end{figure}

The left side of Figure \ref{isosceles type I figure} depicts part of the fingerprint of the point $P$ (in bold) on $X(1,1,3)$ corresponding to the vertex of angle $\dfrac{3\pi}{5}$\,. Since the points corresponding to vertex angles of $\dfrac{\pi}{5}$ (circled) are nonsingular, the sides of the billiards triangulation of $X(1,1,3)$ which are geodesics connecting them to $P$ are not part of the fingerprint of $P$\,.

For a given triangle $T_X$ corresponding to a triangular billiard surface $X$ of genus greater than 1, we can calculate the fingerprint of each element $P$ of a vertex class on $X$\,. We describe this in the following theorem.
\begin{theorem}\label{fingerprint identification theorem}
Suppose $T_X$ unfolds to a surface $X$ of genus greater than 1. Let $T_X$ have vertices $v_1$, $v_2$, and $v_3$\,. Denote the angle measure of $v_i$ by $\alpha_i$\,. Suppose $P \in \pi^{-1}_X(v_1)$\,. Then one of three situations exists:

1) If $v_1$ is the apex of an isosceles triangle, then $P$ has angle set $\{\theta\}$, where $\theta=\alpha_1$\,.

2)  If $v_1$ is not the apex of an isosceles triangle, and if for some $j\in\{2,3\}$ we have that both $\pi_X^{-1}(v_j)$ is singular and $\alpha_j>\dfrac{\pi}{6}$, then $P$ has a Type I fingerprint with angle set $\{\theta\}$, where $\theta=2\alpha_1$\,.

3)If $v_1$ is not the apex of an isosceles triangle, and if for each $j\in\{2,3\}$ either $\alpha_j=\dfrac{\pi}{k}$ for some integer $k>1$ or $\alpha_j < \dfrac{\pi}{6}$, then $P$ has a Type II fingerprint with angle set $\{\theta_1,\theta_2\}$, where $\theta_1=\pi-2\alpha_2$ and $\theta_2=\pi-2\alpha_3$\,.

\end{theorem}
\begin{proof}

1) Suppose that $v_1$ is the apex of an isosceles triangle. Then the shortest geodesics connecting $P$ to a singularity either correspond to the edges of $T_X$ incident on $v_1$ or to the billiard path from $v_1$ to the midpoint of the opposite side and back. See Figure \ref{isosceles type I figure}. The incident edges are longer than the billiard path precisely when $\alpha_2=\alpha_3<\dfrac{\pi}{6}$; in this situation $\alpha_1>\dfrac{2\pi}{3}$ and so $\pi^{-1}(v_1)$ must be singular. Note that when $\alpha_2=\alpha_3=\dfrac{\pi}{6}$ we have that $\pi^{-1}(v_2)$ and $\pi^{-1}(v_3)$ are nonsingular, so the edges and the billiard path cannot both correspond to shortest geodesics connecting $P$ to a singular point. In either case, it is evident that adjacent shortest geodesics form an angle of $\alpha_1$\,.

2) Now suppose that $v_1$ is not the apex of an isosceles triangle, and that for some $j\in\{2,3\}$ we have that both $\pi^{-1}(v_j)$ is singular and $\alpha_j>\dfrac{pi}{6}$\,. Then the edge of $T_X$ with endpoints $v_1$ and $v_j$ is shorter than the shortest billiard path from $v_1$ to itself. Since $v_1$ is not the apex of an isosceles triangle, exactly one of the edges incident on $v_1$ corresponds to all of the shortest geodesics connecting $P$ to singularities. Thus two such geodesics which are adjacent form an angle of $2\alpha_1$\,. See Figure \ref{Type I and II figure}.

3) Finally, suppose that $v_1$ is not the apex of an isosceles triangle, and that for each $j\in\{2,3\}$ either $\alpha_j=\dfrac{\pi}{k}$ for some integer $k>1$ or $\alpha_j < \dfrac{\pi}{6}$\,. It follows that for each side of $T_X$ incident on $v_1$, either the side will not correspond to a geodesic connecting $P$ to a singularity, or else the side will be longer than the shortest billiard path from $v_1$ to itself. Furthermore, we see that $\pi_X^{-1}(v_1)$ must be singular. For, either $\pi_X^{-1}(v_2)$ and $\pi_X^{-1}(v_3)$ are nonsingular, in which case $\pi_X^{-1}(v_1)$ must be singular since $X$ has a singularity, or else $\alpha_2+\alpha_3<\dfrac{2\pi}{3}$; in this case, $\alpha_1>\dfrac{\pi}{3}$ implies that $\pi_X^{-1}(v_1)$ is singular, since if $\alpha=\dfrac{\pi}{2}$ it is impossible to find $\alpha_2$ and $\alpha_3$ satisying the conditions of case 3. The claim is now evident from Figure \ref{X(3,4,5) Type II figure},
which illustrates the fingerprint of the singularity on $X(3,4,5)$\,.
(In the figure, the
geodesics defining
\begin{figure}[ht]
\includegraphics[scale=.35]{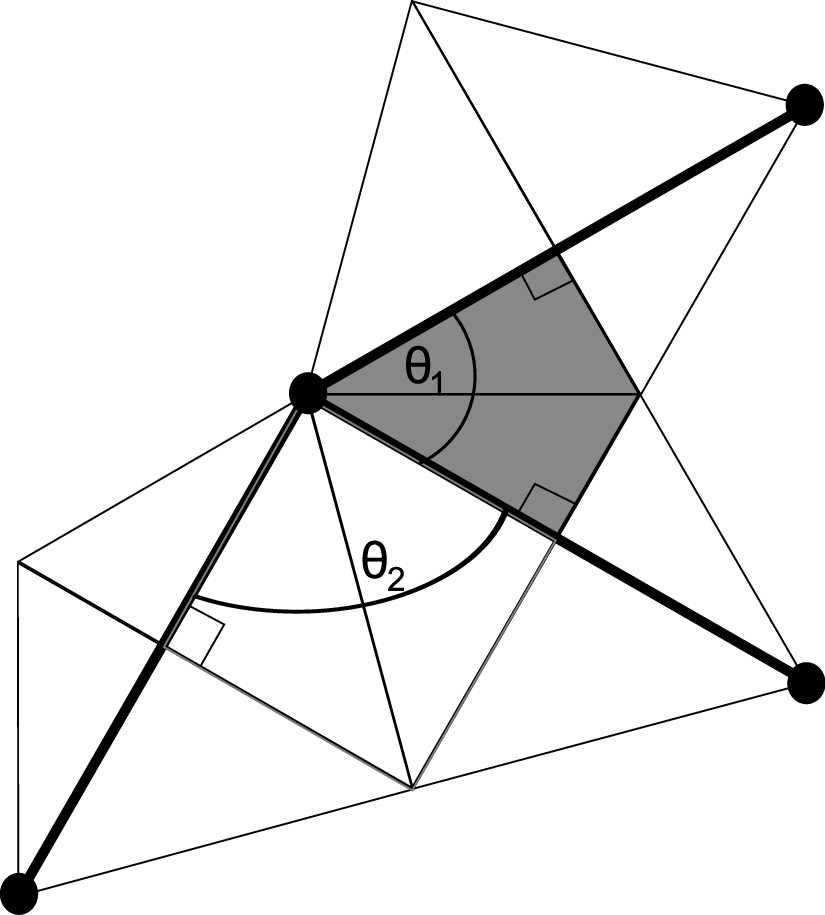}
\caption{Part of a Type II fingerprint on X(3,4,5). Nonsingular vertex points are circled.}
\label{X(3,4,5) Type II figure}
\end{figure}
the fingerprint are the thicker lines, whereas the edges of the
billiards triangulation are the thinner lines.) Let the angle set of the fingerprint of $P$ be
$\{\theta_{1},\theta_{2}\}$\,. Each $\theta_i$ is an interior angle
of a quadrilateral whose other three angles include two right angles
and an angle which has twice the measure of an angle of the
triangular billiard table $T_X$ for $X$\,. Therefore two of the angles
of $T$ have the form
$\dfrac{1}{2}(2\pi-\dfrac{\pi}{2}-\dfrac{\pi}{2}-\theta_i)=\dfrac{\pi-\theta_i}{2}$,
and the third angle is $\dfrac{\theta_{1}+\theta_{2}}{2}$\,.
\end{proof}

\begin{corollary} \label{Type II uniqueness
corollary}
Suppose the billiards triangulation of a
triangular billiard surface $X$ contains a point with a Type II
fingerprint. Then $X$ is uniquely determined by that fingerprint, up
to an action of $O(2,\mathbb{R})$\,. Indeed, if the fingerprint has
angle set $\{\theta_{1},\theta_{2}\}$, then $X$ is the billiards
surface for the triangle of angles
$\dfrac{\theta_{1}+\theta_{2}}{2}$, $\dfrac{\pi-\theta_{1}}{2}$, and
$\dfrac{\pi-\theta_{2}}{2}$\,.
\end{corollary}

The fingerprint is a useful tool for studying translation covers because it is nearly invariant under balanced covers. It seems worth noting that one can define the notion of a fingerprint on any translation surface; this makes particular sense for any polygonal billiard surface. Then the following invariance result will still hold, except that the cone angles may differ by a larger integer factor.

\begin{theorem} \label{fingerprint invariance theorem}
Suppose that $f:X\rightarrow Y$ is a balanced translation
cover, that $P' \in X$ and $P \in Y$ are vertices of billiards
triangulations on their respective surfaces, and that
$f\left(P'\right)=P$\,. Then the fingerprints of $P'$ and $P$ have the same angle sets and shortest geodesic lengths, and their cone angles are either equal or differ only by a factor of two.

\end{theorem}

\begin{proof}
Let $d$ and $d'$ be the lengths of the shortest geodesics which connect $P$ and $P'$, respectively, to a singularity. Since $f$ is a translation cover, the image under $f$ of any geodesic on $X$ is a geodesic of equal or lesser length on $Y$; since $f$ is balanced, any geodesic with singular endpoints on $X$ is mapped by $f$ onto a geodesic with singular endpoints on $Y$\,. Thus, $d\leq d'$\,. Conversely, since $f$ is a translation cover, the preimage of any geodesic of length $d$ with singular endpoints is a union of geodesics of length $d$ with singular endpoints. Thus, $d'\leq d$, and we see that the fingerprints of $P$ and $P'$ have the same geodesic lengths.

Define $B_{\delta}(P)$ to be the set of all points of $Y$ that are of distance less than or equal to $\delta$ from $P$ . Let $B'_{\delta}(P')$ be the connected component of $f^{-1}(B)$ containing $P'$\,. Since $Y$ is a translation surface, and $f$ is a translation cover (possibly ramified above $P$ with index $m$), we can and do choose $\delta$ to be a sufficiently small positive number such that $B_{\delta}$ and $B'_{\delta}(P')$ are simply connected closed metric balls centered at $P$ and $P'$, respectively.   

Consider a pair of adjacent (in the sense of Definition \ref{definition adjacent}) geodesics $e_1$ and $e_2$, each of
length $d$, connecting $P$ to singularities. Label the angle between
them $\theta$\,. Let $l_1$ and $l_2$ be the intersections of $e_1$ and $e_2$ with $B_{\delta}(P)$\,. 

The union of  $l_1$ and $l_2$ with a portion of the
boundary of $B_{\delta}(P)$ bounds a wedge-shaped region $W$ which does not contain in its interior any portion of a shortest geodesic connecting $P$ to a singular point.
Since $f$ is a translation cover, the intersection of the $f^{-1}(W)$ with $B'_{\delta}(P')$ is $m$ copies of $W$\,. Let $W'$ be one of these copies. Then $W'$ is bounded by part of the boundary of $B'$ and two geodesics (say, $l_1'$ and $l_2'$) which extend to be length $d$ geodesics $e_1'$ and $e_2'$ connecting $P'$ to singular points.  Since $f$ is a translation cover, the angle measure between $l_1'$ and $l_2'$ is $\theta$; therefore the angle measure between $e_1'$ and $e_2'$ is $\theta$\,.

Suppose that $e_1'$ and $e_2'$ are not
adjacent. Then there is some geodesic $e_3'$ of length $d$ connecting $P'$ to a singularity. The intersection of this geodesic with $W'$ is a geodesic $l_3'$ in the interior of $W'$\,. But then $f(l_3')$ is a geodesic in the interior of $W$ which extends to the image of $e_3'$ -- a length $d$ geodesic connecting $P$ to a singular point. This contradicts the adjacency of $e_1$ and $e_2$\,. Hence, in fact $e_1'$ and $e_2'$ are adjacent, and and it follows that the fingerprints of $P$ and $P'$ have the same angle sets.

Finally, we claim that $m \leq 2$\,. Let $v$
and $v'$ be the vertices of the triangles $T_Y$ and $T_X$ corresponding
to $P$ and $P'$\,. By Remark \ref{remark on AI combinatorics}, the
cone angle at $P$ is completely determined by $\angle v$\,. But
Theorem \ref{fingerprint identification theorem} tells us that the measure of
$\angle v$ is determined, up to a factor of 2, by the angle set of
the fingerprint of $P$\,. Hence, since the fingerprints of $P$ and
$P'$ have the same angle set, we see that $m \in \{1, 2\}$\,.
\end{proof}

\begin{corollary}\label{fingerprint type invariance corollary}
Fingerprint type is invariant under balanced translation covers.
\end{corollary}

\begin{corollary}\label{Type II no balanced cover corollary}
Any rational triangular billiard surface with a Type II singularity
cannot be a part of any composition of nontrivial balanced covers.
\end{corollary}

\begin{proof}
By  \ref{Type II uniqueness
corollary}. Suppose we have $f:X\rightarrow Y$ a balanced cover
with either $X$ or $Y$ possessing a singularity with a Type II
fingerprint. By Corollary \ref{fingerprint type invariance corollary},
$X$ and $Y$ must both have singularities with Type II fingerprints.
Since a Type II fingerprint identifies the triangular billiards
table of a surface, $X$ and $Y$ must be the same surface.
\end{proof}

As we shall see, the preceding results about fingerprints allow us to quickly classify
all balanced covers in the category of triangular billiards
surfaces. However, to extend our results to unbalanced covers, we shall refine our use of the fingerprint by considering \emph{punctured surfaces}. By deleting a collection of singular points from a surface $X$ to obtain a surface $\tilde{X}$, we will arrive at a set of shortest geodesics connecting a point $P$ in $\tilde{X}$ to the remaining singularities; in general this set will be different than the set of geodesics determining the fingerprint of $P$ on $X$\,. However, since simply deleting singularities of $X$ does not affect the flat metric on open sets not previously containing the deleted points, we can get important information about $X$ by considering the fingerprint of $P$ on $\tilde{X}$\,.

\begin{lemma}\label{fingerprints with punctures lemma}
Let $X$ be a triangular billiard surface with more than one
singular vertex class. Let $\tilde{X}$ be the surface obtained from
$X$ by puncturing either one entire singular vertex class or two
entire singular vertex classes such that neither deleted class
corresponds to an obtuse angle of the triangular billiard table and
such that at least one singular vertex class remains. Let
$\pi_X^{-1}(v_i)$ be a singular vertex class not deleted. Let $P \in
\pi_X^{-1}(v_i)$\,. If $P$ has Type II fingerprint on $\tilde{X}$
with angle set $\{\theta_1,\theta_2\}$, then $X$ arises from
billiards in the triangle with angles $\dfrac{\pi-\theta_1}{2}$,
$\dfrac{\pi-\theta_2}{2}$, and $\dfrac{\theta_1+\theta_2}{2}$\,. If
$P$ has a Type I fingerprint on $\tilde{X}$ with angle set
$\{\theta_1\}$, then the measure of $\angle v_i$ is in
$\{\theta_1,\frac{1}{2}\theta_1\}$\,.
\end{lemma}

\begin{proof}
If none of the punctured points were endpoints of geodesics defining the fingerprint of $P$, then $P$ has the same fingerprint on
$\tilde{X}$ as on $X$, and we are done.

Suppose a singular vertex class has been punctured which contained
endpoints of such geodesics. Then there is a
new ``closest" vertex class to $P$; call it $C$\,. If $C$ does not
contain $P$ then the shortest geodesics connecting $P$ to $C$ are
edges of the billiards triangulation of $X$\,. If $C$ does contain
$P$ then, since a vertex class corresponding to an obtuse angle of
the billiard table must be singular (by Remark \ref{remark on AI
combinatorics}) and we have assumed that no such classes have been
deleted, it follows that the shortest geodesics from $P$ to $C$
correspond to a single reflection in the original dynamical system.
Thus the same reasoning holds as in Theorem \ref{fingerprint identification theorem}.

The only potential difficulty would be if the new ``closest" vertex
class was the one containing $P$, for in that case, since the
shortest geodesics from $P$ to elements of its own class pass
through more than one triangle, we must consider the possibility
that our punctures obstruct these geodesics. However, since the
shortest geodesics are perpendicular to the sides of the triangles
opposite $P$, this is only a problem if the vertex class punctured
is $\pi_X^{-1}(v_j)$ with $\angle v_j=\dfrac{\pi}{2}$\,. But such a
class is nonsingular, so it would not have been punctured.
\end{proof}

\section{All Translation Covers}\label{section classify all covers}

\subsection{Balanced Covers}\label{subsection balanced covers}

 In this subsection we shall prove Theorem
\ref{classification theorem} for balanced covers using the
simple geometric idea of the fingerprint as our main tool (this is Lemma \ref{balanced cover lemma}).
 
Most points on a billiard surface have no more than three shortest geodesics connecting them to singular points.
We define an \emph{exceptional point} on a billiard surface $X$ to be a point which has more than three shortest geodesic paths to singularities. 

\begin{lemma}\label{exceptional is vertex}
Let $X$ be a triangular billiard surface of genus greater than 1. A point $P$ on $X$ is exceptional if and only if $P$ is a vertex of the billiards triangulation of $X$ not corresponding to a right angle.
\end{lemma}
\begin{proof}
If $P$ is a vertex of $\tau$ then Theorem \ref{fingerprint identification theorem} shows that it is exceptional unless it corresponds to a right angle on $T_X$\,. If $P$ does correspond to a right angle then clearly there are two shortest paths to singularities (corresponding to one side of $T_X$ incident on $\pi_X(P)$) unless $T_X$ is isosceles; but then $X=X(1,1,2)$, which has no singularities.

Suppose that $P$ is in the interior of a triangle of $\tau$\,. Then for each vertex of that triangle, $P$ is strictly closer to that vertex than to any other element of its vertex class. Hence, $P$ cannot have more than three shortest paths to singularities. 

Now suppose that $P$ is in the interior of an edge of $\tau$\,. Then $P$ is on an edge of two triangles, so $P$ has two shortest paths of length $L_1$ to members of the opposite vertex class. If $P$ is the midpoint of the edge then it has equidistant unique shortest paths of length $L_2$ to members of the other two vertex classes; but if $L_1=L_2$ then the triangles must be copies of $T(1,1,2)$, and in that case $X$ has no singularities.

\end{proof}

\begin{lemma}\label{balanced preserved exceptional}
A balanced translation cover $f:X \rightarrow Y$ maps the exceptional points of $X$ onto the exceptional points of $Y$\,.
\end{lemma}
\begin{proof}
Balanced covers preserve fingerprints, so this follows from Lemma \ref{exceptional is vertex} and Theorem \ref{fingerprint identification theorem}.
\end{proof}

\begin{lemma}\label{balanced cover lemma}(Theorem 1 restricted to balanced covers)
Let $X$ and $Y$ be triangular billiard surfaces such that the genus
of $X$ is greater than 1. Suppose that $f:X\rightarrow Y$ is a nontrivial balanced
translation cover. Then $f$ is of the form described in Lemma
\ref{right-isosceles lemma}.
\end{lemma}
\begin{proof} Let $f:X \rightarrow Y$ be a balanced cover between triangular billiard surfaces. If $X$ has three exceptional points with distinct fingerprints, then $T_X$ is not isosceles and so by Theorem \ref{fingerprint identification theorem} we know all the angles of $T_X$; the same reasoning holds for $Y$\,. Since balanced covers preserve fingerprints up to cone angle, we see that $X=Y$\,.

If $X$ has only two distinct fingerprints for its exceptional points then so does $Y$, and by Theorem \ref{fingerprint identification theorem} we know that $T_X$ and $T_Y$ are a pair of triangles described in Lemma \ref{right-isosceles lemma}.
\end{proof}

\subsection{Combinatorial Lemmas}\label{subsection combinatorial lemmas}

\begin{lemma}\label{singular inequality lemma}
Let $f:X(a_1,a_2,a_3) \rightarrow X(b_1,b_2,b_3)$ be a degree $n$ translation
cover of triangular billiard surfaces. Then
\begin{equation}
\sum\limits_{a_i\nmid (a_1+a_2+a_3)}a_i \geq n\sum\limits_{b_i\nmid
(b_1+b_2+b_3)}b_i \;.
\end{equation}
\end{lemma}

\begin{proof}
The sum of the cone angles of the singular points of
$X(a_1,a_2,a_3)$ is at least $n$ times the sum of the cone angles of
the singular points of $X(b_1,b_2,b_3)$\,. By Remark \ref{remark on
AI combinatorics}, the result follows.
\end{proof}

\begin{lemma}\label{congruence lemma}
Let $f:X=X(a_1,a_2,a_3) \rightarrow Y=X(b_1,b_2,b_3)$ be a translation
cover of triangular billiard surfaces such that the genus of
$X$ is greater than 1. If \text{
}$a_1+a_2+a_3=b_1+b_2+b_3$ and $f$ is not a composition of covers
from Lemma \ref{right-isosceles lemma}, then $f$ is of degree 1.
\end{lemma}

\begin{proof}
Write $Q_X=Q_Y=Q=a_1+a_2+a_3=b_1+b_2+b_3$\,. Let $n$ be the degree of $f$,
and suppose that $n \geq 2$\,. Lemma \ref{singular inequality lemma}
then gives $\sum\limits_{b_i\nmid Q}b_i \leq \dfrac{Q}{n}$\,. Hence,
since $n \geq 2$, we have

\begin{equation}\label{nonsingular inequality}
\sum\limits_{b_i | Q}b_i \geq \dfrac{Q}{2}.
\end{equation}

 Writing $q_i=\dfrac{Q}{b_i}$, we have the equivalent
expression

\begin{equation}\label{nonsingular inequality with q}
\sum\limits_{b_i | Q}\dfrac{1}{q_i} \geq \dfrac{1}{2}.
\end{equation}

 Note that if $b_i | Q$ then $q_i$ is an integer. Of course, Equation \eqref{nonsingular inequality} is always
satisfied if $T(b_1,b_2,b_3)$ is a right triangle. If
$T(b_1,b_2,b_3)$ is not a right triangle, the equation is rarely
satisfied. Thus we will reduce the problem to two cases.

\begin{case} The triangle $T(b_1,b_2,b_3)$ is not a right triangle. \end{case}
In this case, recalling that $\gcd(b_1,b_2,b_3)=1$, we show that
there are only three possibilities for the $b_{i}$ which satisfy
Equation \eqref{nonsingular inequality}.

If all three $b_i$ divide $Q$ then $Y$ is nonsingular. The only
non-right triangle which unfolds to a nonsingular surface is
$T(1,1,1)$; but since this is also the only triangle with $Q=3$, if
$Y=X(1,1,1)$ then $X=X(1,1,1)$, contradicting our assumption that
$X$ has a singularity.

Hence we can assume for this case that $b_3 \nmid Q$\,. Therefore to
satisfy Equation \ref{nonsingular inequality with q} we seek
integers $q_1,q_2 >2$ such that

\begin{equation}\label{q inequality}
\dfrac{1}{q_1}+\dfrac{1}{q_2}>\dfrac{1}{2}
\end{equation}

 Without loss of generality we assume $q_1 \leq q_2$\,. If
$q_1 \geq 4$, Equation \eqref{q inequality} is impossible. If
$q_1=3$ then Equation \eqref{q inequality} is satisfied if $q_2 \leq
5$\,. Thus the remaining candidates for $Y$ are $X(3,4,5)$ and
$X(3,5,7)$\,. By Lemma \ref{singular inequality lemma}, $X(3,4,5)$
admits at most a degree two cover; by Lemma \ref{lemma cover
combinatorics} the degree two covers satisfying the hypotheses of
the lemma could only be $f:X(2,5,5) \rightarrow X(3,4,5)$ or
$X(1,1,10) \rightarrow X(3,4,5)$\,. However, these maps would have
to be balanced covers, and $X(3,4,5)$ has a singularity with a Type
II fingerprint. Thus by Corollary \ref{Type II no balanced cover
corollary} these maps do not exist. Similarly, the only feasible
cover of $X(3,5,7)$ of degree greater than 1 is $f:X(1,7,7)
\rightarrow X(3,5,7)$; again, this would be a balanced cover, and
$X(3,5,7)$ has a singularity with a Type II fingerprint.

\begin{case} The triangle $T(b_1,b_2,b_3)$ is a right triangle.\end{case}

Assume that $b_1=\dfrac{Q}{2}$\,. If $b_2$ divides $Q$, then by Lemma \ref{right-isosceles lemma}, the surface $Y=X(b_2+b_3,b_2,b_3)$ is translation equivalent to the surface $X(b_3,\frac{b_2}{2},\frac{b_2}{2})$\,. But then this is equivalent to the situation where $Q_X=2Q_Y$; we will address this situation in the proof of Theorem \ref{classification theorem}.

Thus for the remainder of this proof we will assume that neither $b_2$ nor $b_3$ divides $Q$\,.  Here Lemma
\ref{singular inequality lemma} implies that the degree of $f$ is at
most two. The sum of the cone angles of the singularities of $Y$ is
$2\pi(b_2+b_3)$\,. Thus if $n=2$ then the sum of the cone angles of the
singularities of $X$ is $4\pi(b_2+b_3)=2\pi Q=2\pi(a_1+a_2+a_3)$\,. Therefore
$T(a_1,a_2,a_3)$ must be either $T(b_2,b_2,2b_3)$ or
$T(2b_2,b_3,b_3)$\,. Both these possibilities are accounted for by
the covers of Lemma \ref{right-isosceles lemma}.
\end{proof}

\begin{lemma}\label{degree argument lemma}
Let $f: X \rightarrow Y$ be a translation cover of triangular
billiard surfaces. Let $m$ be the smallest integer such that all
singularities of $Y$ have cone angle at least $2m\pi$\,. Suppose
that $\deg f < m$\,. Then for each vertex class $C_i$ on $X$,
$f(C_i)$ consists entirely of singular points or entirely of
nonsingular points.
\end{lemma}

\begin{proof}
Let $m$ be as above and assume that $\deg(f)<m$\,. Suppose for
contradiction that for some $j$, $f(C_j)$ contains singular points
and nonsingular points. Each member of $C_j$ has the same cone
angle, and this cone angle must be at least $2m\pi$, since some of
the members are mapped by a translation cover to a singularity of
cone angle $2m\pi$\,. Thus, for those elements of $C_j$ which are
mapped to nonsingular points, the definition of a ramified cover
requires that $f$ be locally of degree at least $m$, which
contradicts our assumption that $\deg(f)<m$\,. This completes the
proof.
\end{proof}

\setcounter{case}{0}

\subsection{Proof of the Main Theorem}\label{subsection the main theorem}
Now we can prove Theorem \ref{classification theorem}, which for our
ease we restate in the following way.

\setcounter{theorem}{0}
\begin{theorem}
Suppose $f: X \rightarrow Y$ is a translation cover of
triangular billiard surfaces of degree greater than 1, with the genus of $X$ greater than 1. Then $f$ is
of degree 2, and is a composition of one or two of the covers $f_i$
described in Lemma \ref{right-isosceles lemma}.
\end{theorem}

\begin{proof}
Suppose $X:=X(a_1,a_2,a_3)$, $Y:=X(b_1,b_2,b_3)$, and $f:X
\rightarrow Y$ is a translation cover of degree $\deg f>1$\,. Assume
that the genus of $X$ is greater than 1. Write $Q_X:=a_1+a_2+a_3$
and $Q_Y:=b_1+b_2+b_3$\,. Let $v_1,v_2,v_3$ and $w_1,w_2,w_3$ be the
corresponding vertices of $T(a_1,a_2,a_3)$ and $T(b_1,b_2,b_3)$
respectively.  By Corollary \ref{holonomy field equality corollary}, $X$ and $Y$ have the same holonomy
field $k$\,. By Corollary \ref{holonomy determines Q corollary}, we
have $Q_Y \in \{2Q_X,Q_X,\dfrac{1}{2}Q_X\}$\,. If $Q_Y=2Q_X$, then
by Lemma \ref{singular inequality lemma}, we must have
$\sum\limits_{b_i\nmid Q_Y}b_i \leq
\dfrac{Q_X}{2}=\dfrac{Q_Y}{4}$\,. But then we would have
$\sum\limits_{b_i | Q_Y}b_i \geq \dfrac{3}{4}Q_Y$, which is only the
case for the following surfaces with even $Q$-value: $X(1,1,2)$,
$X(1,2,3)$, $X(3,4,5)$\,. Of course, $Q_X \geq 3$, so $Y \neq
X(1,1,2)$\,. If $Y=X(1,2,3)$ then $X=X(1,1,1)$, which is of genus 1,
a contradiction. If $Y=X(3,4,5)$, then $Y$ has a singularity with
cone angle $10\pi$\,. But, no surface $X$ with $Q_X=6$ could have a
cone angle of at least $10\pi$\,.

If $Q_Y=Q_X$, then we are done by Lemma \ref{congruence lemma}.
Thus, appealing to Corollary \ref{holonomy determines Q corollary},
we shall assume for the remainder of the proof that
$Q_X=2Q_Y$ and that $Q_Y$ is odd.

\begin{case}
$Y$ has three singular vertex classes.
\end{case}

In this case, Lemma \ref{singular inequality lemma} implies that $f$
can only be a degree 2 balanced cover. Thus we are done by Lemma
\ref{balanced cover lemma}.

\begin{case} $Y$ has no singular vertex classes. \end{case}
 In this case, since $Q_Y$ is odd, we
must have $Y=X(1,1,1)$\,. There are only two surfaces with a
$Q$-value of 6: they are $X(1,1,4)$ and $X(1,2,3)$, and each of
these surfaces covers $X(1,1,1)$ as described in Lemma
\ref{right-isosceles lemma}. 

\begin{case}$Y$ has one singular vertex class. \end{case}

In this case we have, without loss of generality, $b_1|Q_Y$,
$b_2|Q_Y$, and $b_3 \nmid Q_Y$\,. Since $b_1$ and $b_2$ are divisors
of the odd number $Q_Y:=b_1+b_2+b_3$, $b_3$ must also be odd.
Therefore $\dfrac{b_3}{\gcd(b_3,Q)} \geq 3$\,. The cone angle at
each of the singularities of $Y$ corresponding to $b_3$ is
$\dfrac{b_3}{\gcd(b_3,Q)}2\pi \geq 6\pi$\,.

Lemma \ref{singular inequality lemma} eliminates all possible $Y$
for $\deg f \geq 4$ except $Y=X(3,5,7)$\,. But, again by Lemma
\ref{singular inequality lemma}, the only possible degree 4 covering
surface would be $X(1,1,28)$, and such a cover would have to be
balanced, contradicting Lemma \ref{balanced cover lemma}.

\underline{If $\deg f =2$:} Lemma \ref{degree argument lemma} tells
us that if $\deg f=2$ then for each $j=1,2,3$, we have that
$f(\pi_X^{-1}(v_j)) \cap \pi_Y^{-1}(w_3)$ is either empty or all of
$f(\pi_X^{-1}(v_j))$\,.

Suppose that $Y=X(3,5,7)$\,. Lemma \ref{degree argument lemma}
restricts the possible degree 2 covers to surfaces of the form
$X(14,a_2,a_3)$, where each of $a_2$ and $a_3$ is either a divisor
of $30$ or twice a divisor of $30$\,. The only possibility this
leaves is $X(15,14,1)$\,. But any translation cover $X(15,14,1)
\rightarrow X(3,5,7)$ would have to be balanced, so Lemma
\ref{balanced cover lemma} applies.

Now suppose that $Y \neq X(3,5,7)$\,. Let $C$ be the singular vertex
class of $Y$\,. We must have $\dfrac{b_3}{Q}>\dfrac{1}{2}$, and so
by Remark \ref{remark on AI combinatorics} $C$ must correspond to an
obtuse angle $\theta$ of the billiard table. Let $\tilde{X}$ be the
surface obtained from $X$ by puncturing all singular vertex classes
of $X$ which are not contained in $f^{-1}(C)$\,. Since
$\dfrac{b_3}{Q}
> \dfrac{1}{2}$ and $f$ is degree $2$, the sum of the angles of the
billiard table corresponding to the vertex classes in the
$f$-preimage of $C$ must be obtuse. Thus we can apply Lemma
\ref{fingerprints with punctures lemma} to $\tilde{X}$\,. The
restriction of $f$ to $\tilde{X}$ is balanced. Since $Y$ has only
one singular vertex class, elements of $C$ must have Type II
fingerprints unless $T(b_1,b_2,b_3)$ is isosceles. If the
fingerprints are Type II, then Proposition \ref{Type II uniqueness
corollary} and Lemma \ref{fingerprints with punctures lemma}
demonstrate that $X$ and $Y$ are translation equivalent. So the only
possibility is that the fingerprints are Type I. In that case $Y$ is
an isosceles triangular billiard surface. Let $C'$ be a vertex
class on $X$ that is in $f^{-1}(C)$, and write
$\theta=\dfrac{b_3\pi}{Q}$\,. The billiard table angle that $C'$
corresponds to is either $\theta$ or $\dfrac{\theta}{2}$\,. If the
angle is $\theta$, then $X$ and $Y$ are translation equivalent. If
the angle is $\theta \backslash 2$, then there is another vertex
class on $X$ which is also mapped to $C$\,. But then that vertex
class would also correspond to an angle of $\theta \backslash 2$,
and we would have that $X$ is an isosceles triangular billiards
surface which covers $Y$ as described in Lemma
\ref{right-isosceles lemma}.

\underline{If $\deg(f)=3$:} Then Lemma \ref{singular inequality
lemma} allows only the following possibilities for $Y$: the surfaces
\\[5 pt]
 $ Y_n = \left \{
\begin{array}
         {r@{\quad \quad}l}
     X(3,n,2n-3) & 3 \nmid n \\\\

     X(1,\frac{n}{3},\frac{2n}{3}-1) & 3 | n
     \end{array} \right.   .\\$
\\[5 pt]
\indent Note that $\gcd(2n-3,3n) \in \{1,3\}$\,. First suppose that
$\gcd(2n-3,3n)=1$\,. Then $Q=3n$ (thus $n$ is odd), $3 \nmid n$, and
we have $Y_n=X(3,n,2n-3)$\,. We have that $n \geq 5$ and hence that
$2n-3 \geq 7$\,. On $Y_n$, there is only one singular vertex class
and the cone angle of each singular point is $(2n-3)2\pi$\,. Thus
Lemma \ref{degree argument lemma} applies here. Since $Y_n$ is never
isosceles, its singular point has a Type II fingerprint. Let
$\tilde{X}$ be the surface obtained from $X$ by deleting all
singularities of $X$ which $f$ maps to nonsingular points, and let
$\tilde{f}$ be the restriction of $f$ to $\tilde{X}$\,. By Lemma
\ref{degree argument lemma}, the elements of $X-\tilde{X}$ are the
union of entire vertex classes. Thus a Type II fingerprint on
$\tilde{X}$ will uniquely identify the triangular billiards table
used to generate $X$, by Lemma \ref{fingerprints with punctures
lemma}. Because $\tilde{f}$ is a balanced map, each singular point
of $\tilde{X}$ must have the same Type II fingerprint (on
$\tilde{X}$) as its $\tilde{f}$-image on $Y$\,. But, a Type II
fingerprint uniquely identifies the triangle used to generate the
surface (this works for $\tilde{X}$ as well); hence $X$ and $Y_n$
are the same billiard surface, and Lemma \ref{singular surface
self-cover lemma} says that a triple cover is impossible.

Now suppose that $\gcd(2n-3,3n)=3$\,. Then the cone angle of the
singular point on $Y_n$ is $\frac{2n-3}{3}2\pi$\,. If $n>6$ then
$\dfrac{2n-3}{3}>3$, so that again we can apply Lemma \ref{degree
argument lemma} and Lemma \ref{fingerprints with punctures lemma},
and the same fingerprint argument goes through. The remaining cases
are $n=3,6$\,. We have $Y_3=X(1,1,1)$ and $Y_6=X(1,2,3)$, neither of
which have singularities.

\begin{case}$Y$ has two singular vertex classes. \end{case}

Assume $b_1|Q$ and $b_2,b_3 \nmid Q$\,. Since $Q$ is odd, we have 
$\dfrac{b_1}{Q} \leq \dfrac{1}{3}$, and so Lemma \ref{singular
inequality lemma} implies that $\deg(f) \leq 3$\,. But, if
$\deg(f)=3$, Lemma \ref{singular inequality lemma} also implies that
$f$ is balanced, contradicting the result of Lemma \ref{balanced
cover lemma} that balanced covers are of degree at most 2. Thus
$\deg(f)=2$\,.

Note that $b_2$ and $b_3$ must have the same parity.



\begin{subcase}
Both $b_2$ and $b_3$ are odd.
\end{subcase}

Then $\dfrac{b_i}{\gcd(b_i,Q)} \geq
3$, so by Lemma \ref{degree argument lemma}, each vertex class of
$X$ maps to all singular points or all nonsingular points.

\emph{If one vertex class of $X$ maps to nonsingular points:} Say
the vertex class $C_1$ corresponding to $a_1$ maps to nonsingular
points. Then $a=2b_1$, and $2b_1|2Q$, so $C_1$ is nonsingular, so
$f$ is balanced.

\emph{If two vertex classes of $X$ map to nonsingular points:}
Let them be $C_1$ and $C_2$, corresponding to $a_1$ and $a_2$\,. If
$C_1$ is singular, then by Lemma \ref{degree argument lemma} we have
$a_1=2d$ for some $d|Q$\,. But since $a_3=2(b_2+b_3)$, this would
mean that all the $a_i$ are even, contradicting the fact that
$\gcd(a_1,a_2,a_3)=1$\,.



\begin{subcase}
Both $b_2$ and $b_3$ are even.
\end{subcase}

\emph{If one vertex class of $X$ maps to nonsingular points:}
Let it be $C_1$\,. We have $a_2+a_3=2(b_2+b_3)$, so $a_1$ must be
even. But also $a_2$ and $a_3$ must be even, since
$2|\dfrac{b_i}{\gcd(b_i,Q)}$ and
$\dfrac{b_i}{\gcd(b_i,Q)}|\dfrac{a_j}{\gcd(a_j,Q)}$ for each $i,j
\in \{2,3\}$\,. Again, this is a contradiction.

\emph{If two vertex classes of $X$ map to nonsingular points:}
Let them be $C_1$ and $C_2$\,. We have that $a_3=2(b_2+b_3)$ is
even. If $C_1$ is singular then again we have that $a_1$ (and hence
$a_2$) is even, once more contradicting that
$\gcd(a_1,a_2,a_3)=1$\,. Hence $C_1$ and $C_2$ are nonsingular, and
$f$ is balanced.

\end{proof}

\end{document}